\documentclass[12pt]{article}
\usepackage{amsmath}
\usepackage{amssymb}
\usepackage{latexsym}
\usepackage{amssymb}
\usepackage{epsf}

\newtheorem{theorem}{Theorem}
\newtheorem{lemma}{Lemma}
\newtheorem{proposition}{Proposition}
\newtheorem{convention}{Convention}
\newtheorem{corollary}{Corollary}
\numberwithin{equation}{section}

\begin{document}
\date{}
\author{M.I.Belishev \thanks {Saint-Petersburg Department of
                 the Steklov Mathematical Institute (POMI), 27 Fontanka,
                 St. Petersburg 191023, Russia; belishev@pdmi.ras.ru.
                 Supported by RFBR grants 08-01-00511 and NSh-4210.2010.1.}}

\title{An unitary invariant of semi-bounded operator and its
application to inverse problems} \maketitle
\begin{abstract}
Let $L_0$ be a closed densely defined symmetric semi-bounded
operator with nonzero defect indexes in a separable Hilbert space
${\cal H}$. With $L_0$ we associate a metric space $\Omega_{L_0}$
that is named a {\it wave spectrum} and constructed from
trajectories $\{u(t)\}_{t \geq 0}$ of a dynamical system governed
by the equation $u_{tt}+(L_0)^*u=0$. The wave spectrum is
introduced through a relevant von Neumann operator algebra
associated with the system. Wave spectra of unitary equivalent
operators are isometric.

In inverse problems on {\it unknown} manifolds, one needs to
recover a Riemannian manifold $\Omega$ via dynamical or spectral
boundary data. We show that for a generic class of manifolds,
$\Omega$ is {\it isometric} to the wave spectrum $\Omega_{L_0}$ of
the minimal Laplacian $L_0=-\Delta|_{C^\infty_0(\Omega\backslash
\partial \Omega)}$ acting in ${\cal H}=L_2(\Omega)$, whereas $L_0$
is determined by the inverse data up to unitary equivalence. By
this, one can recover the manifold by the scheme "the data
$\Rightarrow L_0 \Rightarrow \Omega_{L_0} \overset{\rm isom}=
\Omega$".

The wave spectrum is relevant to a wide class of dynamical
systems, which describe the finite speed wave propagation
processes. The paper elucidates the operator background of the
boundary control method (Belishev, 1986) based on relations of
inverse problems to system and control theory.
\end{abstract}
\setcounter{section}{-1}

\section{Introduction}
\subsection{Motivation}
The paper introduces the notion of a {\it wave spectrum} of a
symmetric semi-bounded operator in a Hilbert space. The impact
comes from inverse problems of mathematical physics; the following
is one of the motivating questions.

Let $\Omega$ be a smooth compact Riemannian manifold with the
boundary $\Gamma$, $-\Delta$ the (scalar) Laplace operator,
$L_0=-\Delta|_{C^\infty_0(\Omega \backslash \Gamma)}$ the {\it
minimal Laplacian} in ${\cal H}=L_2(\Omega)$. Assume that we are
given with a unitary copy $\widetilde L_0=U L_0 U^*$ in
$\widetilde {\cal H} =U{\cal H}$ (but $U$ is unknown!). To what
extent does $\widetilde L_0$ determine the manifold $\Omega$?

So, we have no points, boundaries, tensors, etc, whereas the only
thing given is an operator $\widetilde L_0$ in a Hilbert space
$\widetilde {\cal H}$. Provided the operator is unitarily
equivalent to $L_0$, is it possible to extract $\Omega$ from
$\widetilde L_0$? Such a question is an "invariant" version of
various setups of dynamical and spectral inverse problems on
manifolds \cite{BIP97}, \cite{BIP07}.

\subsection{Content}
Substantially, the answer is affirmative: for a generic class of
manifolds, any unitary copy of the minimal Laplacian determines
$\Omega$ up to isometry (Theorem 1). A wave spectrum is a
construction that realizes the determination $\widetilde L_0
\Rightarrow \Omega$ and, thus, solves inverse problems. In more
detail,
\begin{itemize}
\item With a closed densely defined symmetric semi-bounded
operator $L_0$ of nonzero defect indexes in a separable Hilbert
space ${\cal H}$ we associate a metric space $\Omega_{L_0}$ (its
wave spectrum). The space consists of the so-called eikonal
operators ({\it eikonals}), so that $\Omega_{L_0}$ is a subset of
the bounded operators algebra ${\mathfrak B ({\cal H})}$, whereas
the metric on $\Omega_{L_0}$ is $\|\tau-\tau'\|_{{\mathfrak B
({\cal H})}}$.

The eikonals are constructed from the projections on the reachable
sets of an abstract {\it dynamical system with boundary control}
governed by the evolutionary equation $u_{tt}+L_0^*u=0$. More
precisely, they appear in the framework of a von Neumann algebra
${\frak N}_{L_0}$ associated with the system, whereas
$\Omega_{L_0} \subset {\frak N}_{L_0}$ is a set of the so-called
{\it maximal eikonals}. The peculiarity is that this algebra is
endowed with an additional operation that we call a {\it space
extension}.

Since the definition of $\Omega_{L_0}$ is of invariant character,
the spectra $\Omega_{L_0}$ and $\Omega_{\widetilde L_0}$ of the
unitarily equivalent operators $L_0$ and $\widetilde L_0$ turn out
to be isometric (as metric spaces). So, a wave spectrum is a
(hopefully, new) unitary invariant of a symmetric semi-bounded
operator. \item A wide generic class of the so-called {\it simple
manifolds} is introduced\footnote{Roughly speaking, the simplicity
means that the symmetry group of $\Omega$ is trivial.}. The
central Theorem 1 establishes that for a simple $\Omega$, the wave
spectrum of its minimal Laplacian $L_0$ is isometric to $\Omega$.
Hence, any unitary copy $\widetilde L_0$ of $L_0$ determines the
simple $\Omega$ up to isometry by the scheme $\widetilde L_0
\Rightarrow \Omega_{\widetilde L_0}\overset{\rm isom}=\Omega_{L_0}
\overset{\rm isom}= \Omega$. In applications, it is the procedure,
which recovers manifolds by the BC-method \cite{BIP97},
\cite{BIP07}: the concrete inverse data determine a relevant
$\widetilde L_0$, what enables one to realize the scheme. \item We
discuss one more option: once the wave spectrum of the copy
$\widetilde L_0$ is found, the BC-procedure realizes elements of
the space $\widetilde {\cal H}$ as functions on
$\Omega_{\widetilde L_0}$ \footnote{In the BC-method, such an
option is interpreted as {\it visualization of waves}
\cite{BIP07}.}. Thereafter, one can construct a {\it functional
model} $L_0^{\rm mod}$ of the original Laplacian $L_0$, the model
being an operator in ${\cal H}^{\rm
mod}=L_{2,\,\mu}(\Omega_{\widetilde L_0})$ related with $L_0$
through a similarity (gauge transform). Hopefully, this
observation can be driven to a functional model of a class of
symmetric semi-bounded operators. Presumably, this model will be
{\it local}, i.e., satisfying ${\rm supp\,} L_0^{\rm mod}y
\subseteq {\rm supp\,} y$.
\end{itemize}

\subsection{Comments}
The concept of wave spectrum summarizes rich "experimental
material" accumulated in inverse problems in the framework of the
BC-method, and elucidates operator background of the latter. In
fact, for the first time $\Omega_{L_0}$ has appeared in
\cite{BKac89} in connection with the M.Kac problem; its later
version (called a wave model) is presented in \cite{BIP07} (sec.
2.3.4). Owing to its invariant nature, $\Omega_{L_0}$ promises to
be useful for further applications to unsolved inverse problems of
elasticity theory, electrodynamics, graphs, etc.

Actually, a wave spectrum is an attribute not of a single operator
but an algebra with space extension. In the scalar problems on
manifolds, this algebra is commutative, whereas its wave spectrum
is identical to Gelfand's spectrum of the norm-closed subalgebra
generated by eikonals. However, it is not clear whether this fact
is of general character. The algebras that appear in the above
mentioned unsolved problems, are {\it noncommutative} and the
relation between their wave and Jacobson's spectra is not
understood yet.

By the recent trend in the BC-method, to recover unknown manifolds
via boundary inverse data is to find spectra of relevant algebras
determined by the data \cite{BSobolev}. We hope for further
promotion of this approach.

\section{Wave spectrum}
\subsection{Algebra with space extension}
Let ${\cal H}$ be a separable Hilbert space, ${\rm Lat}{\cal H}$
the lattice of its (closed) subspaces; by $P_{\cal A}$ we denote
the (orthogonal) projection onto ${\cal A} \in {\rm Lat}{\cal H}$.
Also, if ${\cal A}$ is a non-closed lineal set, we put $P_{\cal
A}:=P_{{\rm clos\,}{\cal A}}$. By $\mathfrak B ({\cal H})$ the
bounded operator algebra is denoted.

An one-parameter family $E=\{E^t\}_{t\geq 0}$ of the maps $E^t:
{\rm Lat}{\cal H} \to {\rm Lat}{\cal H}$ is said to be a {\it
space extension} if
\begin{enumerate} \item $E^0\,=\, \rm id$ \item
$E^t\{0\}\,=\,\{0\}, \quad t\geq 0$ \item $t\leq t'$ and ${\cal A}
\subseteq {\cal A}'$ imply $E^t {\cal A} \subseteq E^{t'}{\cal A}'$. \end{enumerate}
It is also convenient to regard $E^t$ as an operation, which
extends the projections, and write $E^t P_{\cal A}=P^t_{\cal A}:=P_{E^t{\cal A}}$.
\smallskip

Assume that an extension $E$ is given; let $a \subset {\rm
Lat}{\cal H}$ be a family of subspaces. By $\mathfrak N [E,a]$ we
denote the minimal von Neumann operator algebra\footnote{i.e., a
unital weekly closed self-adjoint subalgebra of $\mathfrak B
({\cal H})$: see \cite{Mur}}, which contains all projections
$\{P_{\cal A}\,|\, {\cal A} \in a\}$ and is closed with respect to
$E$, i.e., $P\in \mathfrak N [E,a]$ implies $P^t =E^t P \in
\mathfrak N [E,a],\,\,\,t>0$. As is easy to see, such an algebra
is well defined. As any von Neumann algebra, $\mathfrak N [E,a]$
is determined by the set ${\rm Proj\,}\mathfrak N [E,a]$ of its
projections, whereas its additional property is $$E^t\, {\rm
Proj\,}\mathfrak N [E,a] \subset {\rm Proj\,}\mathfrak N [E,a]\,,
\qquad t > 0.$$

Fix a projection $P \in {\mathfrak N}[E,a]$; a positive
self-adjoint operator in ${\cal H}$ of the form
$$\tau_P\,:=\,\int_{0}^\infty t\,dP^t \,,$$
where $P^t=E^t P$, is said to be an {\it eikonal} \footnote{the
term is taken from the motivating applications}; the set of
eikonals is denoted by ${\rm Eik\,}\mathfrak N [E,a]$.

Let us say that we deal with the {\it bounded case} if each
eikonal is a bounded operator (and, hence, belongs to $\mathfrak N
[E,a]$) and the set ${\rm Eik\,}\mathfrak N [E,a]$ is bounded in
$\mathfrak B ({\cal H})$:
\begin{equation}\label{C<infty}
\sup\, \{\|\tau\|\,\,|\,\,\tau \in {\rm Eik\,}\mathfrak N
[E,a]\}\,<\,\infty\,.
\end{equation}
Otherwise, the situation is referred to as {\it unbounded case}
(see the end of sec 3.4).

\begin{convention} Unless otherwise specified, we
deal with the bounded case. \end{convention}

Recall that the set of self-adjoint operators is partially
ordered: for $A,B \in \mathfrak B ({\cal H})$ the relation $A \leq
B$ means $(Ax,x) \leq (Bx,x), \,\,x \in {\cal H}$. Any monotonic
bounded sequence $A_1 \leq A_2 \leq \dots , \,\,\,\sup
\|A_j\|<\infty$ converges in the strong operator topology to
$s$-$\lim A_j \in \mathfrak B ({\cal H})$ (see, e.g.,
\cite{BSol}).

Extend the name {\it eikonal} to all elements of the set $s$-${\rm
clos\,} {\rm Eik\,}\mathfrak N [E,a] \subset {\mathfrak N}[E,a]$
and denote the extended set by the same symbol ${\rm
Eik\,}\mathfrak N [E,a]$. An eikonal $\tau$ is said to be {\it
maximal} if $\tau \geq \tau'$ for any eikonal $\tau'$ comparable
with $\tau$. Let $\Omega_{{\mathfrak N}[E,a]} \subset {\rm
Eik\,}\mathfrak N [E,a]$ be the {\it set of maximal eikonals}.
\begin{lemma}\label{Omeganonempty}
The set $\Omega_{{\mathfrak N}[E,a]}$ is nonempty.
\end{lemma}
{\bf Proof}\,\,\,By (\ref{C<infty}), any totally ordered family of
eikonals $\{\tau_\alpha\}$ has an upper bound $s$-$\overline{\lim}
\tau_\alpha$, which is also an eikonal. Hence, the Zorn lemma
implies $\Omega_{{\mathfrak N}[E,a]}\not= \emptyset$.\,\,\,$\square$
\smallskip

The ${\mathfrak B ({\cal H})}$-norm induces the distance $\|\tau-\tau'\|$ in
$\Omega_{{\mathfrak N}[E,a]}$ and makes it into a metric space, which is
determined by the extension $E$ and the initial reserve of
subspaces $a$. The space $\Omega_{{\mathfrak N}[E,a]}$ is the main subject
of our paper.

\subsection{Extension $E_L$}
Here we introduce the space extension associated with a
semi-bounded self-adjoint operator. Without lack of generality, it
is assumed positive definite: let \begin{equation}\label{L}
L=L^*\,=\,\int^\infty_0 \lambda\,dQ_\lambda\,; \qquad (Ly,y)\geq
\varkappa \,\|y\|^2,\,\,\,y \in {\rm Dom\,} L \subset {\cal H}\,,
\end{equation}
where $d Q_\lambda$ is the spectral measure of $L$ and $\varkappa$
is a positive constant. Such an operator governs the evolution of
a dynamical system
\begin{align}
\label{dynsystem_1} & v_{tt} + L v\,=\,h\,, \qquad t>0\\
\label{dynsystem_2} & v|_{t=0}\,=\,v_t|_{t=0}\,=\,0\,,
\end{align}
where $h \in L_2^{\rm loc}\left((0,\infty); {\cal H}\right)$ is a
${\cal H}$-valued function of time ({\it control}). Its finite
energy class solution $v=v^h(t)$ is represented by the Duhamel
formula
\begin{align}
\notag &  v^h(t)\,=\,\int_0^t L^{-\frac{1}{2}} \sin\left[(t-s)
L^{\frac{1}{2}}\right]\,h(s)\,ds =\,\langle{\rm see (\ref{L})} \rangle\\
\label{trajectory}& =\,\int^t_0 ds \int^\infty_0 \frac{\sin {\sqrt
\lambda}(t-s)}{\sqrt \lambda}\,d Q_\lambda\, h(s)\,, \qquad t \geq
0\,
\end{align}
(see, e.g., \cite{BSol}). In system theory, $v^h(\,\cdot\,)$ is
referred to as a {\it trajectory}; $v^h(t)\in {\cal H}$ is a {\it state}
at the moment $t$. In applications, $v^h$ describes a {\it wave}
initiated by a source $h$. Note that in the case of $\varkappa\leq
0$ the problem (\ref{dynsystem_1}), (\ref{dynsystem_2}) is also
well defined but the representation (\ref{trajectory}) is of
slightly more complicated form. Thus, the assumption $\varkappa>0$
is accepted just for simplicity.

Fix a subspace ${\cal A} \subset {\cal H}$; the set
$${\cal V}_{{\cal A}}^{\,t}:=\{v^h(t)\,|\,\,h \in L_2^{\rm loc}\left((0, \infty); {\cal A}\right)\}\,, \qquad t>0$$
of all states produced by ${\cal A}$-valued controls is called {\it
reachable} (at the moment $t$). Reachable sets increase as ${\cal A}$
increases and/or $t$ grows. Indeed, the representation
(\ref{trajectory}) easily implies
\begin{equation}\label{delay_relation}v^{{\cal T}_\xi h}(t)\,=\,\left({\cal T}_\xi
v^h\right)(t), \qquad t \geq 0\,,
\end{equation}
where ${\cal T}_\xi$ is the right shift operator in $L_2^{\rm
loc}\left((0,\infty); {\cal H}\right)$:
$$\left({\cal T}_\xi
g\right)(t):=\begin{cases}0, \quad \qquad 0\leq &t <\xi\\
g(t-\xi), &t \geq \xi\end{cases}$$ with $\xi \geq 0$. For $0<t
\leq t'$ and ${\cal A} \subseteq {\cal A}'$, we have \begin{equation*}
{\cal V}_{{\cal A}}^{\,t} \ni v^h(t)=\left({\cal T}_{t'-t} v^h\right)(t')=\langle
{\rm see\,} (\ref{delay_relation}) \rangle =v^{{\cal T}_{t'-t}h}(t') \in
{\cal V}_{{\cal A}}^{\,t'} \subseteq {\cal V}_{{\cal A}'}^{\,t'}\,,
\end{equation*}
i.e., the inclusion
\begin{equation}\label{inclusion_V}
{\cal V}_{{\cal A}}^{\,t} \subseteq  {\cal V}_{{\cal A}'}^{\,t'}\,, \qquad 0<t \leq t'
\end{equation}
does hold.
\medskip

Define a family $E_L=\{E^t\}_{t \geq 0}$ of the maps $E^t: {\rm
Lat} {\cal H} \to {\rm Lat} {\cal H}$ by
\begin{equation}\label{definition_EL}
E^0 {\cal A}:={\cal A}, \quad E^t {\cal A}:={\rm clos\,} {\cal V}_{{\cal A}}^{\,t} \,, \quad t > 0\,.
\end{equation}
\begin{lemma}\label{E_is_extension}
$E_L$ is a space extension.
\end{lemma}
{\bf Proof}\,\,\,\,The properties 1 and 2 (see section 1.1) easily
follow from the definitions and the obvious relation
${\cal V}_{\{0\}}^{\,t}=\{0\}$; the property 3 for $0<t \leq t'$ is seen
from (\ref{inclusion_V}). Thus, it remains to verify 3 for $t=0$,
i.e., to check that $0=t\leq t'$ and ${\cal A} \subseteq {\cal A}'$ leads to
$E^{0} {\cal A} \subseteq E^{t'} {\cal A}'$ or, the same, that ${\cal A} \subseteq
E^r{\cal A}$ for all $r>0$ and ${\cal A}\not=\{0\}$.
\smallskip

Take a nonzero $y \in {\cal A}$ and consider (\ref{dynsystem_1}),
(\ref{dynsystem_2}) with the control $h_\varepsilon(t)=\varphi_\varepsilon
(t)\,y$, where $$\varphi_\varepsilon (t):=\begin{cases} 0, &t\in [0,
r-2\varepsilon)\\\frac{1}{\varepsilon^2}, &t \in [r-2\varepsilon, r-\varepsilon)\\-
\frac{1}{\varepsilon^2}, &t \in [r-\varepsilon, r)\\0, &t \in [r, \infty)
\end{cases}$$
($\varepsilon >0$ is small); note that
$$\int_0^r\varphi_\varepsilon(t)\,f(t)\,dt \underset{\varepsilon \to 0} \to - f'(r)$$
for smooth $f$'s, i.e., $\varphi_\varepsilon(t)$ converges to
$\delta^\prime(t-r)$ as a distribution. Define
$$\psi_\varepsilon(\lambda):=\int_0^r
\frac{\sin[\sqrt{\lambda}\,(r-t)]}{\sqrt
\lambda}\,\varphi_\varepsilon(t)\,dt\,=\,\frac{2 \cos(\sqrt{\lambda}
\,\varepsilon)-\cos(\sqrt{\lambda} \,2 \varepsilon) - 1}{\varepsilon^2 \lambda}$$ and
note that $\psi_\varepsilon(\lambda) \underset{\varepsilon \to 0} \to 1$
uniformly w.r.t. $\lambda$ in any compact segment $[\varkappa,
N]$.

Therefore, one has
\begin{align*}
&\left\|y-v^{h_\varepsilon}(r)\right\|^2 = \langle{\rm see\,}
(\ref{trajectory})\rangle=\left\|y-\int_0^r dt
\int_{\varkappa}^\infty \frac{\sin[\sqrt{\lambda}\,(r-t)]}{\sqrt
\lambda}\,dQ_\lambda
[\varphi_\varepsilon(t) y]\right\|^2 =\\
& \left\|y - \int_{\varkappa}^\infty
\psi_\varepsilon(\lambda)\,dQ_\lambda y\right\|^2
=\left\|\int_{\varkappa}^\infty
\left[1-\psi_\varepsilon(\lambda)\right]\,dQ_\lambda
y\right\|^2=\\
& \int_{\varkappa}^\infty |1-\psi_\varepsilon(\lambda)|^2\,d \|Q_\lambda
y\|^2\,\underset{\varepsilon \to 0}\to 0
\end{align*}
by the properties of $\psi_\varepsilon$. The order of integration change
is easily justified by the Fubini Theorem.

Thus, $y=\underset{\varepsilon \to 0}\lim\, v^{h_\varepsilon}(r)$,
whereas $v^{h_\varepsilon}(r)\in E^r{\cal A}$ holds. By the
closeness of $E^r{\cal A}$, we get $y \in E^r{\cal A}$. Hence,
${\cal A} \subseteq E^r{\cal A}$. \,\,\,$\square$
\medskip

So, with each positive definite operator $L$ one associates the
certain space extension $E_L$ by (\ref{definition_EL}).

\subsection{Algebras ${\mathfrak N}_{L, {\cal D}}$ and ${\mathfrak N}_{L_0}$}

Return to the system (\ref{dynsystem_1})--(\ref{dynsystem_2}) and
fix a nonzero subspace ${\cal D} \in {\rm Lat} {\cal H}$ that we'll call a
{\it directional subspace}. It determines a class
\begin{equation}\label{classM_D}{\cal M}_{{\cal D}}:=\left\{h \in C^\infty\left([0, \infty);
{\cal D}\right)\,|\,\,{\rm supp\,} h \subset (0, \infty)
\right\}\end{equation} of smooth ${\cal D}$-valued controls vanishing
near $s=0$, and the sets
\begin{align}
\notag  &{\cal U}^t:=\left\{h(t)- v^{h^{\prime \prime}}(t)\,\biggl |\,h
\in {\cal M}_{{\cal D}}\right\}\,=
\langle\,{\rm see\,} (\ref{trajectory})\, \rangle\,=\\
\label{reachable_sets_U} & \biggl\{h(t)-\int_0^t L^{-\frac{1}{2}}
\sin\left[(t-s) L^{\frac{1}{2}}\right]\, h^{\prime
\prime}(s)\,ds\,\biggl |\,\,h \in {\cal M}_{{\cal D}}\biggr\}, \qquad t \geq 0
\end{align}
(here $(\,\cdot\,)^\prime :=\frac{d}{ds}$), which we also call
{\it reachable} \footnote{The meaning of this definition and term
is clarified later on, when we deal with systems with boundary
control in sec 2.2}. As can be easily derived from
(\ref{delay_relation}), the sets ${\cal U}^t$ increase as $t$ grows.
\smallskip

Now, take $$a_{L,{\cal D}}\,:=\,\{{\rm clos\,} {\cal U}^t\}_{t \geq 0}\,\subset
\,{\rm Lat}{\cal H}$$ in capacity of the initial family of subspaces
(see sec 1.1). The pair $L, {\cal D}$ determines the algebra with
extension
$${\mathfrak N}_{L, {\cal D}}:={\mathfrak N}[E_L, a_{L,{\cal D}}]$$ and its eikonals  ${\rm
Eik\,}{\mathfrak N}_{L, {\cal D}}$. Assuming that (\ref{C<infty})
holds for ${\mathfrak N}_{L, {\cal D}}$, the set of maximal
eikonals $$\Omega_{L, {\cal D}}\,:=\,\Omega_{{\mathfrak N}_{L,
{\cal D}}}$$ is well defined and called a {\it wave spectrum} of
the pair $L, {\cal D}$. The operator
\begin{equation}\label{boundary eikonal}
\tau^\partial:=\int_{0}^\infty t\,dP_{{\cal U}^t}
\end{equation} is
said to be a {\it boundary eikonal}. The set
\begin{equation}\label{boundary of spectrum}\partial \Omega_{L,
{\cal D}}:=\left\{\tau \in \Omega_{L, {\cal D}}\,|\, \tau \geq
\tau^\partial\right\}\,\subset\, \Omega_{L, {\cal
D}}\end{equation} is a {\it boundary} of the wave spectrum.

Thus, each pair $L, {\cal D}$ determines an algebra with space
extension ${\mathfrak N}_{L, {\cal D}}$ and all corresponding
attributes.
\bigskip

Let $L_0$ be a closed densely defined symmetric semi-bounded
operator with {\it nonzero} defect indexes $n_{\pm} = n \leq
\infty$. As is easy to see, such an operator is necessarily
unbounded. For the sake of simplicity, it is assumed positive
definite: $(L_0 y, y)\geq \varkappa \|y\|^2,\,\,y \in {\rm Dom}
L_0$ with $\varkappa>0$. Let $L$ be the extension of $L_0$ by
Friedrichs, so that $L=L^* \geq \varkappa\,\mathbb I$ and
\begin{equation}\label{L0subLsubL0*} L_0
\subset L \subset L_0^*\end{equation}
holds \cite{BSol}. Also,
note that $1 \leq {\rm dim\,}{\rm Ker\,}L_0^*= n \leq \infty$.
Taking
$${\cal D}\,:=\,{\rm Ker\,}L_0^*$$ as a directional subspace, we can
constitute the pair $L, {\rm Ker\,}L_0^*$, which determines the
algebra
$${\mathfrak N}_{L_0}:={\mathfrak N}_{L, {\rm Ker\,}L_0^*}$$ and its eikonals ${\rm
Eik\,}{\mathfrak N}_{L_0}$. If (\ref{C<infty}) holds, the set of maximal
eikonals
$$\Omega_{L_0}:=\Omega_{L, {\rm Ker\,}L_0^*}$$ is well defined and
referred to as a {\it wave spectrum} of the operator $L_0$; its
subset $\partial \Omega_{L_0}$ is a {\it boundary} of the wave
spectrum.

So, with every $L_0$ of the above-mentioned class, one associates
the algebra ${\mathfrak N}_{L_0}$ and the algebra eikonals ${\rm
Eik\,}{\mathfrak N}_{L_0}$. If the latter set is bounded, the
operator $L_0$ possesses the wave spectrum $\Omega_{L_0}\not=
\emptyset$ \footnote{However, there are examples in applications,
in which $\Omega_{L_0}$ consists of a single point.}.

\section{DSBC}
\subsection{Green system}
Dynamical systems with boundary control (DSBC) are defined in the
next sections 2.2 and 2.3; here we introduce a basic ingredient of
the definition. The ingredient is a collection $\{ {\cal H},{\cal
G}; A,\Gamma_0, \Gamma_1\}$ of the separable Hilbert spaces ${\cal
H}$ and ${\cal G}$, and the densely defined operators $A:{\cal H}
\to {\cal H}, \, \Gamma_k:{\cal H} \to {\cal G}\,\,\,(k=0,1)$
connected through the Green formula
$$(Au,v)_{\cal H} - (u,Av)_{\cal H} =
(\Gamma_0 u,\Gamma_1 v)_{\cal G} - (\Gamma_1 u, \Gamma_0 v)_{\cal
G}$$ (see \cite{Kot}). Such a collection is said to be a {\it
Green system}; ${\cal G}$ and $\Gamma_k$ are referred to as a {\it
boundary values space} (BVS) and the {\it boundary operators}
respectively. In the applications, which we deal with later on,
the following is also provided:
\begin{enumerate}
\item ${\rm Dom\,} \Gamma_k \supseteq {\rm Dom\,} A$ holds, whereas $A$ is such
that the restriction $$A|_{{\rm Ker\,} \Gamma_0 \cap {\rm Ker\,}
\Gamma_1}=:L_0$$ is a densely defined symmetric po\-si\-tive
definite operator with nonzero defect indexes and $\overline
A=L_0^*$ is valid ('bar' is the operator closure) \item the
restriction $$A|_{{\rm Ker\,} \Gamma_0} =: L$$ coincides with the
Friedrichs extension of $L_0$, so that we have
\begin{equation}\label{L0subLsubL0*subA}
L_0 \subset L \subset L_0^*=\overline A\,,
\end{equation}
whereas $L^{-1}$ is bounded and defined on ${\cal H}$\item for the
subspaces ${\cal A}:={\rm Ker\,} A$ and ${\cal D}:={\rm Ker\,} L_0^*$, the relations
\begin{equation}\label{density} {\rm clos\,} {\cal A} = {\cal D}\,, \qquad {\rm clos\,}
\Gamma_0 {\cal A}={\cal G}\end{equation} are valid.
\end{enumerate}
These properties are in consent with the Green system theory by
V.A.Ryzhov, which puts them as the basic axioms. Note, that there
are a few versions of such an axiomatics but the one proposed in
\cite{Ryzh} is most relevant for applications to the forward and
inverse {\it multidimensional} problems of mathematical physics.
\begin{convention}
In what follows we deal with the Green systems, which satisfy the
conditions 1--3.
\end{convention}
As is shown in \cite{Ryzh}, the axioms provide the following:
\begin{itemize}
\item the map $\Pi:=\left(\Gamma_1 L^{-1}\right)^*: {\cal G} \to {\cal H}$ is
bounded, whereas ${\rm Ran\,} \Pi$ is dense in ${\cal D}$ \item the subspace
${\cal A}$ admits the characterization
\begin{equation}\label{charact_A}{\cal A}=\{y \in {\rm Dom\,} A\,|\,\,\Pi
\Gamma_0 y = y\}
\end{equation}
\item since $L$ is the extension of $L_0$ by Friedrichs, we have
\begin{equation}\label{L0 via L and D}
{\rm Dom\,} L_0 = L^{-1}[{\cal H} \ominus {\cal D}]\,, \quad  L_0\,=\,L|_{L^{-1}[{\cal H}
\ominus {\cal D}]}
\end{equation}
that easily follows from the definition of such an extension (see
\cite{BSol}).
\end{itemize}
\smallskip

{\bf Example}\,\,\,\,Let $\Omega$ be a $C^\infty$-smooth compact
Riemannian manifold with the boundary $\Gamma$, $\Delta$ the
(scalar) Laplace operator in ${\cal H}:=L_2(\Omega)$, $\nu$ the outward
normal on $\Gamma$, ${\cal G}:=L_2(\Gamma)$. Denote
\footnote{$H^k(\,\dots\,)$ are the Sobolev classes;
$H^2_0(\Omega)=\{y \in H^2(\Omega)\,|\,y=|\nabla y|=0\,\, {\rm
on}\,\Gamma\}$.}
\begin{equation*}A=-\Delta|_{H^2(\Omega)}\,, \qquad \Gamma_0
:=(\,\cdot\,)|_\Gamma\,, \qquad \Gamma_1
:=\frac{\partial}{\partial \nu}(\,\cdot\,)|_\Gamma\,,
\end{equation*}
so that $\Gamma_{0, 1}$ are the trace operators. The collection
$\{ {\cal H},{\cal G}; A,\Gamma_0, \Gamma_1\}$ is a Green system,
whereas other operators, which appear in the framework of Ryzhov's
axiomatics, are the following:
$$L_0=-\Delta|_{H^2_0(\Omega)}$$ is the minimal Laplacian;
$$L=-\Delta|_{H^2(\Omega)\cap H^1_0(\Omega)}$$ is the self-adjoint
Dirichlet Laplacian; $$L_0^*=-\Delta |_{\{y \in {\cal H}\,|\,\Delta y
\in {\cal H}\}}$$ is the maximal Laplacian;
$${\cal A} = \{y \in H^2(\Omega)\,|\,\,\Delta y =0\}$$ is the set of
harmonic functions of the class $H^2(\Omega)$; $${\cal D}=\{y \in
{\cal H}\,|\,\,\Delta y=0\}$$ is the subspace of all harmonic functions
in $L_2(\Omega)$; $\Pi: {\cal G} \to {\cal H}$ is the harmonic continuation
operator (the Diriclet problem solver): $$\Pi \varphi=u: \qquad
\Delta u=0\,\,\,{\rm
in}\,\,\,\Omega,\,\,\,u|_{\Gamma}=\varphi\,.$$

\subsection{Evolutionary DSBC}
Note in advance that the goal of the sections 2.2, 2.3 is not to
obtain results but motivate the above introduced objects and
notions. That is why the presentation ignores certain technical
details.
\medskip

{\bf Definition}\,\,\,\,The Green system determines an
evolutionary {\it dynamical system with boundary control} of the
form
\begin{align}
\label{0.1} & u_{tt}+A\,u=0 &&{\rm in}\,\, {\cal H}, \quad 0<t<\infty   \\
\label{0.2} & u|_{t=0}=u_t|_{t=0}=0 &&{\rm in}\,\, {\cal H} \\
\label{0.3} & \Gamma_0 u=f(t) &&{\rm in}\,\,{\cal G}, \quad 0\leq
t < \infty,
\end{align}
where $f$ is a {\it boundary control} (${\cal G}$-valued function
of time), $u=u^f(t)$ is the solution ({\it wave}).

Assign $f$ to a class ${\cal F}_+$ if it belongs to $C^\infty
\left([0, \infty); {\cal G}\right)$, takes the values in $\Gamma_0
{\rm Dom\,} A$, and vanishes near $t=0$, i.e., satisfies ${\rm
supp\,} f \subset (0, \infty)$. Also, note that $f \in {\cal F}_+$
implies $\Pi \left(f(\,\cdot\,)\right) \in {\cal M}_{{\cal D}}$,
where ${\cal D}= {\rm Ker\,} L_0^*$ and ${\cal M}_{{\cal D}}$ is
defined by (\ref{classM_D}).
\begin{lemma}
For $f \in {\cal F}_+$, the classical solution $u^f$ to problem
(\ref{0.1})--(\ref{0.3}) is represented in the form
\begin{equation}\label{u^f=h-int}
u^f(t)\,=\,h(t)-\int_0^t L^{-\frac{1}{2}} \sin\left[(t-s)
L^{\frac{1}{2}}\right]\, h^{\prime \prime}(s)\,ds\,, \qquad t \geq
0
\end{equation}
with $h:=\Pi \left(f(\,\cdot\,)\right) \in {\cal M}_{{\cal D}}$.
\end{lemma}
{\bf Proof}\,\,\,\,Introducing a new unknown
$w=w^f(t):=u^f(t)-\Pi\left(f(t)\right)$ and taking into account
(\ref{charact_A}), we easily get the system
\begin{align*}
& w_{tt}+A w=-\Pi\left(f_{tt}(t)\right) &&{\rm in}\,\, {\cal H}, \quad 0<t<\infty   \\
& w|_{t=0}=w_t|_{t=0}=0 &&{\rm in}\,\, {\cal H} \\
& \Gamma_0 w = 0 &&{\rm in}\,\,{\cal G}, \quad 0\leq t <\infty\,.
\end{align*}
With regard ro the definition of the operator $L$ (see the axiom
2), this problem can be rewritten in the form
\begin{align*}
& w_{tt}+L w\,=\,- h_{tt} &&{\rm in}\,\, {\cal H}, \quad 0<t<\infty   \\
& w|_{t=0}=w_t|_{t=0}=0 &&{\rm in}\,\, {\cal H}
\end{align*}
and then solved by the Duhamel formula
\begin{equation*}
w^f(t)\,=-\int_0^t L^{-\frac{1}{2}} \sin\left[(t-s)
L^{\frac{1}{2}}\right]\, h^{\prime \prime}(s)\,ds\,.
\end{equation*}
Returning back to $u^f=w^f + \Pi f$, we arrive at
(\ref{u^f=h-int}). \,\,\,\,$\square$
\smallskip

The sets
\begin{align}
\notag & {\cal U}^t_+\,:=\,\{u^f(t)\,|\,\, f \in {\cal F}_+\}=
\langle\,{\rm see\,} (\ref{u^f=h-int})\, \rangle\,=\\
\label{reachable_from_boundary} &\biggl\{h(t)-\int_0^t
L^{-\frac{1}{2}} \sin\left[(t-s) L^{\frac{1}{2}}\right]\,
h^{\prime \prime}(s)\,ds\,\biggl |\,\,\,\,h=\Pi
f(\,\cdot\,),\,\,\,f \in {\cal F}_+ \biggr\}\,, \quad t \geq 0
\end{align}
are said to be {\it reachable from boundary}. In the mean time,
the Green system, which governs the DSBC, determines the certain
pair $L, {\cal D}$, which in turn determines the family $\{{\cal U}^t\}$ by
(\ref{reachable_sets_U}). Comparing the definitions, we easily
conclude that the inclusion ${\cal U}^t_+ \subset {\cal U}^t$ holds. Moreover,
the density properties (\ref{density}) enable one to derive
\begin{equation}\label{clUxipartial=clUxi}{\rm clos\,} {\cal U}^t_+ = {\rm clos\,} {\cal U}^t\,, \qquad t
\geq 0\end{equation} and it is the relation, which inspires the
definition (\ref{reachable_sets_U}) and motivates the use of the
term "reachable set" for ${\cal U}^t$ in the general case, where
neither the boundary value space nor the boundary operators are
defined.
\bigskip

{\bf Illustration}\,\,\,\,Consider the Example at the end of sec
2.1. The DSBC (\ref{0.1})--(\ref{0.3}) associated with the
Riemannian manifold is governed by the wave equation and is of the
form
\begin{align}
\label{0.1Rim} & u_{tt}-\Delta u=0 &&{\rm in}\,\,\,\Omega \times (0, \infty) \\
\label{0.2Rim} & u|_{t=0}=u_t|_{t=0}=0 &&{\rm in}\,\,\,\Omega \\
\label{0.3Rim} & u|_\Gamma =f(t) &&{\rm for}\,\,\,0\leq t < \infty
\end{align}
with a boundary control $f \in {\cal F} := L_2^{\rm loc}\left((0,
\infty); L_2(\Gamma) \right)$; the solution $u=u^f(x,t)$ describes
a wave, which is initiated by boundary sources and propagates from
the boundary into the manifold with the speed 1. For $f \in
{\cal F}_+:=C^\infty \left([0, \infty); C^\infty(\Gamma) \right)$
provided ${\rm supp\,} f \subset (0, \infty)$, the solution $u^f$ is
classical. By the finiteness of the wave propagation speed, at a
moment $t$ the waves fill the near-boundary subdomain
$$\Omega^t[\Gamma]:=\{x \in \Omega\,|\,\,{\rm dist\,} (x, \Gamma) <t\}\,.$$
Correspondingly, the reachable sets ${\cal U}^t_+$ increase as $t$ grows
and the relation
\begin{equation}\label{Uxi_sub_Hxi}
{\cal U}^t_+ \,\subset {\cal H}^t\,, \qquad t \geq 0
\end{equation}
holds, where ${\cal H}^t := {\rm clos\,} \{y \in {\cal H}\,|\,\,{\rm supp\,} y \subset
\Omega^t[\Gamma]\}$.
\bigskip

So, if the pair $L, {\cal D}$ (or the operator $L_0$) appears in
the framework of a Green system, then $\{{\cal U}^t\}$ introduced
by the general definition (\ref{reachable_sets_U}) may be imagine
as the sets of waves produced by boundary controls. The question
arises, what is the meaning of the corresponding wave spectrum
$\Omega_{L, {\cal D}}$ ($=\Omega_{L_0}$)? In a sense, it is the
question, which this paper is written for. The answer (postponed
till section 3) is that $\Omega_{L_0}$ is a wave guide body, in
which such waves propagate.
\bigskip

{\bf Controllability}\,\,\,\,Return to the abstract DSBC
(\ref{0.1})--(\ref{0.3}) and define its certain property. The
definition is premised with the following observation. Since the
class of controls ${\cal F}_+$ satisfies $\frac{d^2}{dt^2}{\cal F}_+={\cal F}_+$,
the reachable sets (\ref{reachable_from_boundary}) satisfy $A
{\cal U}^t_+ ={\cal U}^t_+$. Indeed, taking $f \in {\cal F}_+$ we have
\begin{equation}\label{relations}
A u^f(t)=\langle\,{\rm see\,} (\ref{0.1})\,\rangle = - u^f_{tt}(t)
= u^{-f^{\prime \prime}}(t) \in {\cal U}^t_+\,
\end{equation}
and, by the same relations, $$u^f(t)\,=\,A u^g(t)$$ with
$g=-(\int_0^t)^2 f \in {\cal F}_+$. Hence, the sets ${\cal U}^t_+$ reduce the
operator $A$, so that its parts $A|_{{\cal U}^t_+}$ are well defined.
\smallskip

The DSBC (\ref{0.1})--(\ref{0.3}) is said to be {\it controllable}
from boundary for the time $t=T$ if the (operator) closure of the
part $A|_{{\cal U}^T_+}$ coincides with $\overline A$, i.e., the
relation
\begin{equation}\label{controllability_def_1}
{\rm clos\,} \left\{\{u^f(T), A u^f(T)\}\,|\,\,f \in
{\cal F}_+\right\}\,=\,{\rm clos\,} {\rm graph\,} A
\end{equation}
(the closure in ${\cal H} \times {\cal H}$) is valid, where $${\rm graph\,}
A:=\left\{\{y, A y\}\,|\,\,y \in {\rm Dom\,} A \right\}\,.$$
Controllability means two things. First, since $A$ is densely
defined in ${\cal H}$, (\ref{controllability_def_1}) implies
$${\rm clos\,} {\cal U}^t_+ \,=\,{\cal H}\,, \qquad t \geq T\,,$$ i.e.,
for large times the reachable sets become rich enough: dense in
${\cal H}$. Second, the "wave part" $A|_{{\cal U}^T_+}$ of the operator $A$,
which governs the evolution of the system, represents the operator
in substantial: it coincides with $A$ up to closure.

In applications to problems in bounded domains, such a property
"ever holds" (for large enough times $T$). In particular, the
system (\ref{0.1Rim})--(\ref{0.3Rim}) is controllable for any
$T>\underset{x \in \Omega} \max\, {\rm dist\,}(x, \Gamma)$ (see
\cite{BIP97}, \cite{BIP07}).
\medskip

Let us represent the property (\ref{controllability_def_1}) in the
form available for what follows. At first, with regard to
(\ref{L0subLsubL0*subA}) and (\ref{relations}), it can be written
as
\begin{equation}\label{controllability_def_2}
{\rm clos\,} \left\{\{u^f(T), u^{-f^{\prime \prime}}(T)\}\,\bigg|\,\,f
\in {\cal F}_+\right\}\,={\rm graph\,} L_0^*\,.
\end{equation}
Further, for each $t \geq 0$, introduce a {\it control operator}
("input $\to$ state" map) $W^t: {\cal F} \to {\cal H}, \,\,\,{\rm Dom\,} W^t = {\cal F}_+$,
$$W^t f\,:=\, u^f(t)\,.$$ In terms of this map, (\ref{controllability_def_2})
takes the form
\begin{equation}\label{controllability_def_3}
{\rm clos\,} \left\{\{W^T f, W^T (-f^{\prime \prime})\}\,\bigg|\,\,f \in
{\cal F}_+\right\}\,={\rm graph\,} L_0^*\,.
\end{equation}
The control operator can be regarded as an operator from a Hilbert
space ${\cal F}^t:=L_2\left((0, t); {\cal G}\right)$ to the space ${\cal H}$ with
${\rm Dom\,} W^t={\cal F}^t_+:=\{f|_{[0, t]}\,\,|\,\, f \in {\cal F}_+\}$. As such,
it can be represented in the form of the {\it polar decomposition}
$$W^t\,=\,U^t\,|W^t|\,,$$ where
$|W^t|:=\left((W^t)^*W^t\right)^{1 \over 2}$ and $U^t$ is an
isometry from ${\rm clos\,} {\rm Ran\,} |W^t| \subset {\cal F}^t$ onto ${\rm clos\,} {\rm Ran\,} W^t
\subset {\cal H}$ (see, e.g., \cite{BSol}). For $t=T$, one has ${\rm clos\,}
{\rm Ran\,} W^T = {\cal H}$, so that $U^T$ is a unitary operator from the
(sub)space $\widetilde {\cal H}:={\rm clos\,} {\rm Ran\,} |W^T|$ onto ${\cal H}$.
Correspondingly, the operator $$\widetilde
L_0^*\,:=\,\left(U^T\right)^*L_0^*\, U^T $$ acting in $\widetilde
{\cal H}$ turns out to be unitarily equivalent to $L_0^*$. As result,
(\ref{controllability_def_3}) can be written in the final form
\begin{equation}\label{controllability_def_4}
{\rm clos\,} \left\{\left\{|W^T| f,\, |W^T| (-f^{\prime
\prime})\right\}\,\bigg|\,\,f \in {\cal F}_+\right\}\,={\rm graph\,} \widetilde
L_0^*\,.
\end{equation}
As a consequence, we get
\begin{proposition}\label{|WT|_determines_L0*}
If the DSBC (\ref{0.1})--(\ref{0.3}) is controllable for the time
$T$ then the operator $|W^T|$ determines the operator $L_0^*$ up
to unitary equivalence.
\end{proposition}
\bigskip

{\bf Response operator}\,\,\,\,In the DSBC
(\ref{0.1})--(\ref{0.3}), an "input$\to$output" correspondence is
described by the {\it response operator} $R: {\cal F} \to {\cal F}, \,\,\,{\rm Dom\,}
R = {\cal F}_+$,
$$\left(Rf\right)(t)\,:=\,\Gamma_1 \left(u^f(t)\right)\,, \qquad t
\geq 0\,.$$ Also, the reduced operators
$$ R^t f\,:=\,\left(Rf\right) |_{[0, t]}$$ are in use and play
the role of the data in dynamical inverse problems.

The key fact of the BC-method is that the operator $R^{2t}$
determines the operator $(W^t)^*W^t$ through a simple and explicit
relation: see \cite{BIP97}, \cite{DSBC}, \cite{BIP07}. Hence,
$R^{2t}$ determines the modulus $|W^t|$. Combining this fact with
the Proposition \ref{|WT|_determines_L0*}, we arrive at
\begin{proposition}\label{R2T_determines_L0*}
If the DSBC (\ref{0.1})--(\ref{0.3}) is controllable for the time
$T$ then its response operator $R^{2T}$ determines the operator
$L_0^*$ (and, hence, the operator $L_0=L_0^{**}$ and its
Friedrichs extension $L$) up to unitary equivalence.
\end{proposition}
\medskip

As illustration, the response operator of the DSBC
(\ref{0.1Rim})--(\ref{0.3Rim}) is $$R^{2T}: f\,\mapsto
\frac{\partial u^f}{\partial \nu}\bigg|_{\Gamma \times
[0,2T]}\,.$$ By the aforesaid, given for a fixed $T>\underset{x
\in \Omega} \max\, {\rm dist\,}(x, \Gamma)$ this operator determines the
operator $L_0$ up to a unitary equivalence.

\subsection{Stationary DSBC}
Our presentation follows the paper \cite{Ryzh}. The basic object
is the Green system $\{ {\cal H},{\cal G}; A,\Gamma_0, \Gamma_1\}$
and the associated operators $L_0, L$ (see sec 2.1).
\medskip

{\bf Definition}\,\,\,\,Along with the evolutionary DSBC, one
associates with the Green system the problem
\begin{align}
\label{0.1Stat} & \left(A - z \mathbb I\right) u=0 &&{\rm in}\,\,
{\cal H},\,\,\, z \in \mathbb C  \\
\label{0.2Stat} & \Gamma_0 u=\varphi  &&{\rm in}\,\,{\cal G}
\end{align}
that is referred to as a {\it stationary} DSBC.  For $\varphi \in
\Gamma_0 {\rm Dom\,} A$ and $z \in {\mathbb C\,}\backslash {\,\rm spec\,}
L$, such a problem has a unique solution $u=u^\varphi_z$, which is
a ${\rm Dom\,} A\,$-valued function of $z$.
\medskip

{\bf Weyl function}\,\,\,\,The "input $\to$ output" correspondence
in the system (\ref{0.1Stat})--(\ref{0.2Stat}) is realized by an
operator-valued function $$M(z) \varphi \,:=\,\Gamma_1
u^\varphi_z\,, \qquad z \notin {\,\rm spec\,} L$$ that is called
{\it Weyl function} and plays the role of the data in frequency
domain inverse problems.

The following fact proven in \cite{Ryzh} is of crucial value.
Recall that a symmetric operator in ${\cal H}$ is said to be {\it
completely non-selfadjoint} if there is no subspace in ${\cal H}$, in
which the operator induces a self-adjoint part.
\begin{proposition}\label{M_determines_L0}
If the Green system, which determines the DSBC
(\ref{0.1Stat})--(\ref{0.2Stat}), is such that the operator $L_0$
is completely non-selfadjoint, then the Weyl function of the DSBC
determines the operator $L_0$ up to unitary equivalence.
\end{proposition}
\bigskip

{\bf Illustration}\,\,\,\,Consider the Example at the end of sec
2.1. The DSBC (\ref{0.1Stat})--(\ref{0.2Stat}) associated with the
Riemannian manifold is
\begin{align}
\label{0.1StatRim} & \left(A + z \right) u=0 \qquad {\rm in}\,\,\Omega\\
\label{0.2StatRim} & u|_\Gamma\,=\,\varphi\,,
\end{align}
where $A=-\Delta|_{H^2(\Omega)}$. The operator
$L_0=-\Delta|_{H^2_0(\Omega)}$ is completely non-selfadjoint.
Indeed, otherwise there exists a subspace ${\cal K} \subset {\cal H}$ such
that the operator $L_0^{{\cal K}}:=-\Delta|_{{\cal K} \cap H^2_0(\Omega)}\not=
\mathbb O$ is self-adjoint in ${\cal K}$. In the mean time, $L_0^{{\cal K}}$
is a part of $L$, the latter being a self-adjoint operator with
the discrete spectrum. Hence, ${\rm spec}\,L^{{\cal K}}_0$ is also
purely discrete; each of its eigenfunctions satisfies $-\Delta
\phi=\lambda \phi$ in $\Omega$ and belongs to $H^2_0(\Omega)$. The
latter implies $\phi = \frac{\partial \phi}{\partial \nu}=0$ on
$\Gamma$, which leads to $\phi \equiv 0$ by the well-known
E.Landis uniqueness theorem for solutions to the Cauchy problem
for elliptic equations. Hence, $L^{{\cal K}}_0=\mathbb O$ in
contradiction to the assumptions.

The Weyl function of the system  is $$M(z) \varphi\,=\,
\frac{\partial u^\varphi_z}{\partial \nu}\bigg|_\Gamma \qquad (z
\not \in {\rm spec\,} L)\,.$$ By the aforesaid, the function $M$
determines the operator $L_0$ up to a unitary equivalence.
\smallskip

Besides the Weyl function, there is one more kind of inverse data
associated with to the DSBC
(\ref{0.1StatRim})--(\ref{0.2StatRim}). Let
$\{\lambda_k\}_{k=1}^\infty : \,\,\,0<\lambda_1<\lambda_2\leq
\lambda_3 \leq \dots \to \infty$ be the spectrum of the Dirichlet
Laplacian $L$, $\{\phi_k\}_{k=1}^\infty :
\,\,\,L\phi_k=\lambda_k\phi_k$ its eigenbasis in ${\cal H}$
normalized by $(\phi_k, \phi_l)=\delta_{kl}$. The set of pairs
$$\Sigma_\Omega\,:=\,\left\{\lambda_k ;\,\frac{\partial
\phi_k}{\partial \nu}\bigg|_\Gamma\right\}_{k=1}^\infty$$ is
called the (Dirichlet) {\it spectral data} of the manifold
$\Omega$. The well-known fact is that these data determine the
Weyl function and vice versa (see, e.g., \cite{Ryzh}). Hence,
$\Sigma_\Omega$ determines the minimal Laplacian $L_0$ up to
unitary equivalence. However, such a detrmination can be realized
not through $M$ but in more explicit way.

Namely, let $U: {\cal H} \to \widetilde {\cal H} := \emph{l}_2$,
\begin{equation*}
U y\,=\,\widetilde y:=\{(y,\phi_k)\}_{k=1}^\infty\end{equation*}
be the Fourier transform that diagonalizes $L$:
\begin{equation}\label{diag_L}
\widetilde L\,:=\,U L U^*={\rm diag\,}\{\lambda_1,\,\lambda_2\,,\,
\dots\}\,.\end{equation} For any harmonic function $a \in {\cal A}$, its
Fourier coefficients are
\begin{equation*}
(a,\phi_k)\,=\,-\,\frac{1}{\lambda_k} \int_\Gamma
a\,\frac{\partial \phi_k}{\partial \nu}\,d\Gamma
\end{equation*}
that can be verified by integration by parts. With regard to the
latter, the spectral data $\Sigma_\Omega$ determine the image
$\widetilde {\cal A}:=U {\cal A} \subset \widetilde {\cal H}$ and
its closure $\widetilde {\cal D}=U {\cal D}={\rm clos\,}
\widetilde {\cal A}$, so that the determination
\begin{equation*}
\Sigma_\Omega\,\Rightarrow \widetilde L\,, \widetilde {\cal D}
\end{equation*}
occurs. In the mean time, (\ref{L0 via L and D}) implies
\begin{equation}\label{tildeL0}
\widetilde L_0\,=\,U^* L_0 U\,=\, \widetilde L |_{{\widetilde
L}^{-1}\left[\widetilde {\cal H} \ominus \widetilde {\cal D}\right]}
\end{equation}
by isometry of $U$. Thus, $\widetilde L_0$ is a unitary copy of
$L_0$ constructed via the spectral data.

\section{Applications}

\subsection{Inverse problems}
In inverse problems (IP) for DSBC associated with manifolds, one
needs to recover the manifold via its boundary inverse data
\footnote{In concrete applications (acoustics, geophysics,
electrodynamics, etc), these data formalize the measurements
implemented at the boundary.}. Namely,
\smallskip

\noindent {\bf IP 1:}\,\,\,given for a fixed $T>\underset{x \in
\Omega} \max\, {\rm dist\,}(x, \Gamma)$ the response operator
$R^{2T}$ of the system (\ref{0.1Rim})--(\ref{0.3Rim}), to recover
the manifold $\Omega$
\smallskip

\noindent {\bf IP 2:}\,\,\,given the Weyl function $M$ of the
system (\ref{0.1StatRim})--(\ref{0.2StatRim}), to recover the
manifold $\Omega$
\smallskip

\noindent {\bf IP 3:}\,\,\,given the spectral data
$\Sigma_\Omega$, to recover the manifold $\Omega$.
\smallskip

\noindent The problems are called {\it time-domain}, {\it
frequency-domain}, and {\it spectral} respectively.

Setting the goal to determine an unknown manifold from its
boundary inverse data, we have to keep in mind the evident
nonuiqueness of such a determination: all {\it isometric}
manifolds with the mutual boundary have the same data. Therefore,
the only reasonable understanding of "to recover" is to construct
a manifold, which possesses the prescribed data \cite{BIP07}.

As we saw, the common feature of the problems IP 1--3 is that
their data determine the minimal Laplacian $L_0$ up to unitary
equivalence. By this, each kind of data determines the wave
spectrum $\Omega_{L_0}$ up to isometry. As will be shown, for a
wide class of manifolds the relation $\Omega_{L_0} \overset{\rm
isom}= \Omega$ holds. Hence, for such manifolds, for solving  the
IPs it suffices to extract a unitary copy $\widetilde L_0$ from
the data, find its wave spectrum $\Omega_{\widetilde L_0}
\overset{\rm isom}=\Omega_{L_0}$, and thus to get an isometric
copy of $\Omega$. It is the program for the rest of the paper.

\subsection{Simple manifolds}
Recall that we deal with a compact smooth Riemannian manifold
$\Omega$ with the boundary $\Gamma$; $\rm vol$ is the volume in
$\Omega$. Also, recall some definitions.

For a subset $A \subset \Omega$, denote by
$$\Omega^r[A]\,:=\,\{x \in \Omega\,|\,\,{\rm dist\,}(x, A)<r\}$$ the
metric neighborhood of $A$ of radius $r>0$ and put
$\Omega^0[A]:=A$. Note that whatever $A$ be, its neighborhood is
an open set with the zero volume boundary:
\begin{equation}\label{Federer}{\rm vol\,\,}\partial
\Omega^r[A]\,=\,0\,,\qquad r>0\end{equation} \cite{Fed}. By
$A^\flat$ we denote the set of its interior points: $x \in
A^\flat$ if there is an $\varepsilon>0$ such that
$\Omega^\varepsilon [\{x\}]\subset A$. For a system $\alpha
\subset 2^\Omega$, we define $\alpha^\flat:=\{A^\flat\,|\,\,A \in
\alpha\}$.

Let $Y$ be a set, $\Xi \subset 2^Y$ a system of subsets. The
system $\Xi$ is said to be an {\it algebra\,\,\,} if
\begin{itemize}\item $Y \in \Xi\,\,\,$ \item $A, B \in \Xi\,\,\,$
implies $\,\,\,Y\backslash A, \,A \cap B \in \Xi$ (and hence
$\,\,\,\emptyset,\,A \cup B \in \Xi$).\end{itemize} For a family
$\alpha \subset 2^Y$, by $\Xi[\alpha]$ we denote the algebra
generated by this family, i.e., the {\it minimal algebra} that
contains $\alpha$. As is known, $\Xi[\alpha]$ consists of the sets
of the form $\cup_{n=1}^N \cap_{m=1}^M A_{nm}$, where $A_{mn}$ or
$Y\backslash A_{mn}$ belong to $\alpha$ (see, e.g., \cite{BSol}).
\medskip

Return to the manifold. The following is an universal process that
associates with $\Omega$ a certain system $\alpha_\Gamma \subset
2^\Omega$ and that we refer to as a {\it Procedure $1$}. The
process is a consequent repetition of the same operation $\sigma$
that acts as follows: for a given family $\alpha \subset
2^\Omega$,
\begin{enumerate}
\item constitute the algebra $\Xi[\alpha]$ and go to the system
$\Xi^\flat[\alpha]$\item construct the system of the neighborhoods
$$\sigma[\alpha]:=\{\Omega^t[A]\,|\,\,A \in \Xi^\flat[\alpha],\,\,t
\geq 0\}$$ that is the product of the operation $\sigma$.
\end{enumerate} Such an operation is of the following important
feature. We say a set $A \subset \Omega$ to be {\it regular} (and
write $A \in {\cal R} \subset 2^\Omega$) if ${\rm vol\,\,}A>0$ and ${\rm
vol\,\,}\partial A=0$ holds. Note that ${\cal R}$ is an algebra. The
feature is that by the definition of the minimal algebra and
property (\ref{Federer}), the inclusion $\alpha \subset {\cal R}$
implies $\sigma [\alpha] \subset {\cal R}$. Also, note that the passage
$\Xi[\alpha] \to \Xi^\flat[\alpha]$ removes the zero volume sets,
which can appear in the algebra $\Xi[\alpha]$.

Now, we describe
\smallskip

\noindent  {\bf Procedure $\bf 1$}:

\noindent{\bf Step 1}\,\,\,\,Take the family of boundary
neighborhoods $\gamma:=\{\Omega^t[\Gamma]\}_{t\geq 0}\subset
2^\Omega$ and construct the system $\sigma[\gamma]$

\noindent{\bf Step 2}\,\,\,\,Construct $\sigma^2[\gamma]\,:=\,
\sigma[\sigma[\gamma]]$

\noindent{\bf Step 3}\,\,\,\,Construct $\sigma^3[\gamma]\,:=\,
\sigma[\sigma[\sigma[\gamma]]]$

\noindent \dots \quad \dots \quad \dots

\noindent{\bf Final Step}\,\,\,\, Constitute the system
$$\alpha_\Gamma\,:=\,\bigcup \limits_{j=1}^\infty \sigma^j[\gamma]$$
that is the end product of the Procedure $1$. As is easy to see,
the constructed system consists of regular sets, and is determined
by the metric in $\Omega$ and the "shape" of its boundary
$\Gamma$.
\medskip

A system $\alpha \subset 2^\Omega$ is said to be a {\it net} if
for any point $x \in \Omega$ there exists a sequence
$\{\omega_j\}_{j=1}^\infty \subset \alpha$ such that $\,\,{\rm
vol\,} \omega_j>0,\,\,\,\omega_1 \supset \omega_2 \supset \,
\dots\,\,$, ${\rm diam\,}\omega_j \to 0$, and $x \in \bigcap
\limits_{j \geq 1}\omega_j$.
\medskip

We say the manifold $\Omega$ to be {\it simple}, if the system
$\alpha_\Gamma$ is a net. The following is some comments on this
definition.

The evident obstacle for a manifold to be simple is its
symmetries\footnote{Presumably, any compact manifold with trivial
symmetry group is simple, but it is a conjecture. Note that for
noncompact manifolds this is not true.}. For a ball $\Omega=\{x
\in {\mathbb R}^n\,|\,\,|x|\leq 1\}$, the system $\alpha_\Gamma$
consists of the sets $\eta \times S^{n-1}$, where $\eta \subset
[0,1]$ is a sum of positive measure segments. Surely, such a
system is not a net in the ball. A plane triangle is simple iff
its legs are pair-wise nonequal. Sufficient and easily checkable
conditions on the shape of $\Omega \subset {\mathbb R}^n$, which
provide the simplicity, are proposed in \cite{BKac89}. They are
also available for Riemannian manifolds and show that simplicity
is a generic property: it can be reached by arbitrarily small
smooth variations of the boundary $\Gamma$.
\medskip

For an $A \subset \Omega$, define a distant function
$$d_A(x)\,:= \,{\rm dist\,}(x, A)\,, \qquad x \in \Omega\,.$$
By the well-known properties of the distance on a metric space,
distant functions are continuous: $d_A \in C(\Omega)$. The latter
space is a Banach algebra (with the $\sup$-norm). Let
$C_\Gamma(\Omega)$ be the (closed) subalgebra in $C(\Omega)$
generated by the family $\{d_A\,|\,\,A \in \alpha_\Gamma \}$. The
following is the result, which in fact inspires the notion of
simplicity.
\begin{lemma}\label{CGamma=COmega}
If the manifold $\Omega$ is simple then the equality
\begin{equation}\label{C_Gamma(Omega)=C(Omega)}
C_\Gamma(\Omega)\,=\,C(\Omega)
\end{equation}
holds.
\end{lemma}
{\bf Proof}\,\,\,\,For $x, x' \in \Omega,\,\,x\not=x'$, choose
$\omega \in \alpha_\Gamma$ such that $x \in \omega$ and $x' \not
\in \overline \omega$, what is possible since the system
$\alpha_\Gamma$ is a net. One has $0=d_\omega(x)$ and
$d_\omega(x')>0$, so that $C_\Gamma(\Omega)$ distinguishes points
of $\Omega$. Therefore, by the Stone theorem, the sum
$C_\Gamma(\Omega) \vee \{\rm constants\}$ coincides with
$C(\Omega)$. In the mean time, assuming $\{\rm constants\} \not
\subset C_\Gamma(\Omega)$, we have $C(\Omega)=C_\Gamma(\Omega)
\overset{.}+\{\rm constants\}$ and conclude that
$C_\Gamma(\Omega)$ is a maximal ideal in $C(\Omega)$. By the
latter, there exists a (unique) point $x_0 \in \Omega$ such that
all functions of $C_\Gamma(\Omega)$ vanish at $x_0$ (see, e.g.,
\cite{Mur}). Evidently, it is not the case since $\alpha_\Gamma$
is a net. Hence, we arrive at
(\ref{C_Gamma(Omega)=C(Omega)}).\,\,\,$\square$

\subsection{Solving IPs}
Here we prove the basic
\begin{theorem}\label{Theorem_1}
Let $\Omega$ be a simple manifold, $L_0=-\Delta|_{H^2_0(\Omega)}$
the minimal Laplacian, $\Omega_{L_0}$ its wave spectrum. There
exists an isometry (of metric spaces) $i$ that maps $\Omega_{L_0}$
onto $\Omega$, the relation $i(\partial \Omega_{L_0})=\Gamma$
being valid.
\end{theorem}
{\bf Proof}\,\,consists of the parts I--III.
\medskip

{\bf I.\,\,\,System $\hat \alpha_\Gamma$}

We say a system of subspaces $\hat\Xi \subset {\rm Lat\,}{\cal H}$ to be
an {\it algebra\,\,\,} if \begin{itemize} \item ${\cal H} \in \hat
\Xi\,\,\,$ \item ${\cal A}, {\cal B} \in \hat \Xi\,\,\,$ implies $\,\,\,{\cal H}
\ominus {\cal A}, \,{\cal A} \cap {\cal B} \in \hat \Xi$ (and hence
$\,\,\,\{0\},\,{\cal A} \vee {\cal B} \in \hat \Xi$).\end{itemize} For a
family $a \subset {\rm Lat\,}{\cal H}$, by $\hat \Xi[a]$ we denote the
algebra generated by this family, i.e., the {\it minimal algebra}
that contains $a$. As is known, $\hat \Xi[a]$ consists of the
subspaces of the form $\cup_{n=1}^N \vee_{m=1}^M {\cal A}_{nm}$, where
${\cal A}_{mn}$ or ${\cal H} \ominus {\cal A}_{mn}$ belong to $a$.

Below we present an universal process that associates with
$\Omega$ a certain system of subspaces $\hat\alpha_\Gamma \subset
{\rm Lat\,}{\cal H}$ and that we refer to as a {\it Procedure
$\hat 1$}. The process is a consequent repetition of the same
operation $\hat \sigma$ that acts as follows: for a given family
$a \subset {\rm Lat\,}{\cal H}$,
\begin{enumerate}\item  constitute the algebra
$\hat \Xi[a]$ \item construct the family
$$\hat \sigma[a]:=\{E^r{\cal A}\,|\,\,{\cal A} \in \hat \Xi[a],\,\,r
\geq 0\}\,\subset {\rm Lat\,}{\cal H}\,,$$ where $\{E^r\}_{r \geq
0}=E_L$ is the space extension determined by the Dirichlet
Laplacian $L \supset L_0$.
\end{enumerate} The latter family is the product of the operation $\hat\sigma$.
Such an operation is of the following important feature.

A subspace ${\cal A} \subset {\cal H}$ is called {\it regular} (we
write ${\cal A} \in \hat {\cal R} \subset {\rm Lat\,}{\cal H})$ if
$${\cal A} =\{y \in {\cal H}\,|\,\,{\rm supp\,} y \subset
\overline A\}\,=:\, {\cal H} A$$ with an $A \in {\cal R}$; so, a
regular subspace consists of functions supported on a regular set.
Now, we invoke a fundamental property of the DSBC
(\ref{0.1Rim})--(\ref{0.3Rim}) known as a {\it local
controllability} \footnote{see \cite{BIP07}, sec 2.2.3, eqn
(2.21). This property is based upon the fundamental
Holmgren-John-Tataru theorem on uniqueness of continuation of
solutions to the wave equation (\ref{0.1Rim}) through a
noncharacteristic surface: see \cite{BIP97} for detail.}, by which
\begin{itemize}
\item the embedding (\ref{Uxi_sub_Hxi}) is dense, i.e., in our
current notation, we have ${\rm clos\,} {\cal U}^t_+ = {\cal H} \Omega^t[\Gamma]$,
which shows that all ${\rm clos\,} {\cal U}^t_+$ are regular subspaces. As a
consequence, by (\ref{clUxipartial=clUxi}) we conclude that the
subspaces ${\rm clos\,} {\cal U}^t$ are also regular:
\begin{equation}\label{Holmgren1}{\rm clos\,} {\cal U}^t\, = \,{\cal H}
\Omega^t[\Gamma]\,\subset \hat {\cal R}, \qquad t
>0\end{equation}\item
for any regular subspace ${\cal H} A$, the relation
\begin{equation}\label{Holmgren2}
E^t {\cal H} A\,=\,{\cal H}\Omega^t[A]\,, \qquad t \geq 0 \end{equation}
holds. \end{itemize} The above-announced feature of the operation
$\hat \sigma$ is the following. By the definition of the minimal
algebra and property (\ref{Holmgren2}), if $\hat \sigma$ is
applied to a family $a$ of regular subspaces, then the result
$\hat \sigma[a]$ also consists of regular subspaces, i.e., $\hat
\alpha \subset \hat {\cal R}$ implies $\hat \sigma[\hat \alpha] \subset
\hat {\cal R}$.

Now, we describe
\smallskip

\noindent{\bf Procedure $\bf \hat 1$}:

\noindent{\bf Step 1}\,\,\,\,Take the family of subspaces $\hat
\gamma:=\{{\rm clos\,} {\cal U}^t\}_{t\geq 0}\subset \hat {\cal R}$ corresponding to
the boundary neighbourhoods $\Omega^t[\Gamma]$ (see
(\ref{Holmgren1})) and construct the family $\hat \sigma[\hat
\gamma]$

\noindent{\bf Step 2}\,\,\,\,Construct
$\hat\sigma^2[\hat\gamma]\,:=\,
\hat\sigma[\hat\sigma[\hat\gamma]]$

\noindent{\bf Step 3}\,\,\,\,Construct
$\hat\sigma^3[\hat\gamma]\,:=\,
\hat\sigma[\hat\sigma[\hat\sigma[\hat\gamma]]]$

\noindent \dots \quad \dots \quad \dots

\noindent{\bf Final Step}\,\,\,\, Constitute the family
$$\hat\alpha_\Gamma\,:=\,\bigcup \limits_{j=1}^\infty \hat\sigma^j[\hat\gamma]$$
that is the end product of the Procedure $\hat 1$. As is easy to
see, the constructed family consists of regular subspaces: $\hat
\alpha_\Gamma \subset \hat {\cal R}$.
\medskip

The evident duality of the Procedures $1$ and $\hat 1$, which is
intentionally emphasized by the notation, easily leads to the
bijection
\begin{equation}\label{bijection} {\cal R} \supset \alpha_\Gamma \ni
A\,\,\longleftrightarrow\,\,{\cal H} A \in \hat \alpha_\Gamma \subset
\hat {\cal R}\,.\end{equation}
\bigskip

{\bf II.\,\,\,Eikonals}\,\,\, Recall that for a linear set ${\cal A}
\subset {\cal H}$, by $P_{\cal A}$ we denote the projection in ${\cal H}$ onto
${\rm clos\,} {\cal A}$.

By the bijection (\ref{bijection}), each projection $P_{\cal A}$
with ${\cal A} \in \hat \alpha_\Gamma$ is of the form $P_{\cal
A}=P_{{\cal H} A}$ with $A \in \alpha_\Gamma$. Hence, as an
operator in ${\cal H} = L_2(\Omega)$, $P_{\cal A}$ multiplies
functions by the indicator $\chi_A (\,\cdot\,)$ of the set $A$
(cuts off functions on $A$), what we write as $P_{\cal A}=\chi_A$.
By (\ref{Holmgren2}), the projection $P^t_{\cal A}=E^tP_{\cal A}$
cuts off functions on $\Omega^t[A]$, i.e., $P^t_{\cal
A}=\chi_{\Omega^t[A]}$. As result, the operator
$$\tau_{P_{\cal A}}\,=\,\int_{0}^\infty t\,dP^t_{\cal A}\,=\,\int_{0}^\infty
t\,d\chi_{\Omega^t[A]}$$ multiplies functions by the distant
function: \begin{equation}\label{tau=d_A}\left(\tau_{P_{\cal A}}
y\right)(x)\,=\,d_A(x)\,y(x)\,, \qquad x \in
\Omega\,.\end{equation}

Let ${\mathfrak L}_\infty(\Omega)$ be the algebra of
$L_\infty$-multipliers, which is a von Neumann subalgebra of
${\mathfrak B}({\cal H})$. By ${\mathfrak C}_\Gamma(\Omega)$ and
${\mathfrak C}(\Omega)$ we denote the subalgebras of
$C_\Gamma(\Omega)$- and $C(\Omega)$- multipliers respectively,
both of them being closed w.r.t. the operator norm; so, we have
${\mathfrak C}_\Gamma(\Omega) \subset {\mathfrak C}(\Omega)
\subset{\mathfrak L}_\infty(\Omega) \subset{\mathfrak B}({\cal
H})$. These subalgebras are commutative; also, as is well known,
${\mathfrak C}(\Omega)$ is weakly dense in ${\mathfrak
L}_\infty(\Omega)$ \cite{Mur}. Recall that ${\mathfrak N}_{L_0}$
is defined in sec 1.3, whereas $L_0$ which we are dealing with, is
the minimal Laplacian on $\Omega$.
\begin{lemma}\label{NL_0=L_infty}
If the manifold $\Omega$ is simple then one has ${\mathfrak
N}_{L_0}={\mathfrak L}_\infty(\Omega)$.
\end{lemma}
{\bf Proof} (sketch)\,\,\, By construction of $\hat
\alpha_\Gamma$, all projections $P_{\cal A}$ with ${\cal A} \in \hat
\alpha_\Gamma$ belong to the algebra ${\mathfrak N}_{L_0}$. Whence, the
eikonals of the form (\ref{tau=d_A}) belong to ${\mathfrak N}_{L_0}$, so
that the embedding ${\mathfrak C}_\Gamma(\Omega) \subset
{\mathfrak N}_{L_0}$ occurs.

By simplicity of $\Omega$, one has ${\mathfrak
C}_\Gamma(\Omega)={\mathfrak C}(\Omega)$, what is just a form of
writing the assertion of Lemma \ref{CGamma=COmega}. So, we have
${\mathfrak C}(\Omega) \subset {\mathfrak N}_{L_0}$. The latter, by the
above mentioned w-closeness of ${\mathfrak C}(\Omega)$ in
${\mathfrak L}_\infty(\Omega)$, implies ${\mathfrak
L}_\infty(\Omega)\subset {\mathfrak N}_{L_0}$.

All projections belonging to ${\mathfrak L}_\infty(\Omega)$ are of
the form $P=\chi_A$ with a Borel $A$ \cite{Mur}. Using the
relevant generalization of the relation (\ref{Holmgren2}) on Borel
positive volume sets, one can show that the algebra ${\mathfrak
L}_\infty(\Omega)$ is invariant w.r.t. the extension $E_L$. Since
${\mathfrak N}_{L_0}$ is a minimal $E_L$-invariant algebra, the embedding
${\mathfrak L}_\infty(\Omega)\subset {\mathfrak N}_{L_0}$ yields
${\mathfrak L}_\infty(\Omega)= {\mathfrak N}_{L_0}$. \,\,\,\,$\square$
\begin{corollary}
The set of eikonals ${\rm Eik\,}{\mathfrak N}_{L_0}$ consists of the
operators of the form (\ref{tau=d_A}) with the Borel $A$'s.
\end{corollary}
\bigskip

{\bf III.\,\,\,Wave spectrum}\,\,\,We omit the proof of the
following simple fact: the eikonal $\tau_{P_{\cal A}} =d_A$ is maximal
iff $A$ is a single point set, i.e., $A=\{x_0\}$ for a $x_0 \in
\Omega$. Such an eikonal is denoted by $\tau_{x_0}$, so that we
have $$\Omega_{L_0}\,=\,\{\tau_{x_0}\,|\,\,x_0 \in \Omega\}\,$$ by
the definition of a wave spectrum, whereas a map
$$i: \,\Omega_{L_0}\, \ni \tau_{x_0} \mapsto x_0 \in\, \Omega$$ is
a bijection. Recall that the distance between eikonals
$\tau',\tau'' \in \Omega_{L_0}$ is $\|\tau'-\tau''\|$ and show
that $i$ is an isometry.

Take $x',x'' \in \Omega$; the corresponding eikonals $\tau_{x'},
\tau_{x''} \in \Omega_{L_0}$ multiply functions by $d_{x'}$ and $
d_{x''}$ respectively. For any $x \in \Omega$, by the triangular
inequality one has $|d_{x'}(x) - d_{x''}(x)| \leq {\rm dist\,} (x',
x'')$; hence for $y \in {\cal H}=L_2(\Omega)$ one easily has
$\|\left[\tau_{x'} - \tau_{x''}\right] y\| \leq {\rm dist\,} (x',
x'')\|y\|$ that implies
\begin{equation}\label{||tau'-tau''||<}
\|\tau_{x'} - \tau_{x''}\| \leq {\rm dist\,} (x', x'')\,.
\end{equation}
Choose $\omega_j \subset \Omega$ such that $\omega_1 \supset
\omega_2 \supset \dots\,, \, {\rm diam\,} \omega_j \to 0$, and $x'
= \cap_{j \geq 1}\, \omega_j$; then put
$y_j=\|\chi_j\|^{-1}\chi_j$, where $\chi_j$ is the indicator of
$\omega_j$. As is easy to see, $\|\tau_{x'} y_j\| \to 0$ and
$\|\tau_{x''} y_j\| \to {\rm dist\,} (x', x'')$ holds as $j \to \infty$
that implies
\begin{equation}\label{||tau'-tau''||to dist}
\|\left[\tau_{x'} - \tau_{x''}\right] y_j\| \to {\rm dist\,} (x', x'')\,,
\qquad \|y_j\|=1.
\end{equation}
Comparing (\ref{||tau'-tau''||<}) with (\ref{||tau'-tau''||to
dist}), we arrive at $\|\tau_{x'} - \tau_{x''}\|= {\rm dist\,} (x', x'')$
and conclude that the map $\tau_{x_0} \mapsto x_0$ is an isometry.
\medskip

Recall that the boundary of the wave spectrum is introduced in sec
1.3 by (\ref{boundary eikonal}) and (\ref{boundary of spectrum}).
In our case, in accordance with (\ref{Holmgren1}) the projections
$P_{{\cal U}^t}$ cut off functions on the near-boundary subdomains
$\Omega^t[\Gamma]$ and, correspondingly, the boundary eikonal
$\tau^\partial=\int_0^\infty t\,d P_{{\cal U}^t}$ multiplies functions
by $d_\Gamma:={\rm dist\,} (\,\cdot\,, \Gamma)$. As is evident, the the
relation $d_{x_0} \geq d_\Gamma$ holds in $\Omega$ iff $x_0 \in
\Gamma$. Equivalently, the eikonal $\tau_{x_0} \in \Omega_{L_0}$
satisfies $\tau_{x_0}\geq \tau^\partial$ iff $x_0 \in \Gamma$.
This means that the isometry $i$ maps the boundary $\partial
\Omega_{L_0}$ of the wave spectrum onto the boundary $\Gamma$ of
the manifold.

Theorem 1 is proved. \,\,\,\,$\square$
\bigskip

Regarding non-simple manifolds, note the following. If the
 symmetry group of $\Omega$ is nontrivial then, presumably,
$\Omega_{L_0}$ is isometric to the properly metricized set of the
group orbits. Such a conjecture is supported by the following
easily verifiable examples:
\begin{itemize}
\item for a ball $\Omega=\{x \in {\mathbb R}^n\,|\,\,|x| \leq
r\}$, the spectrum $\Omega_{L_0}$ is isometric to the segment
$[0,r] \subset \mathbb R$, whereas its boundary $\partial
\Omega_{L_0}$ is identical to the endpoint $\{0\}$\item for an
ellipses $\Omega=\{(x, y)\in {\mathbb R}^2\,|\,\,\frac{x^2}{a^2} +
\frac{y^2}{b^2} \leq 1\}$, $\Omega_{L_0}$ is isometric to its
quarter $\Omega \cap \{(x, y)\,|\,\,x \geq 0,\, y \geq 0 \}$,
whereas $\partial \Omega_{L_0}\overset{\rm isom}=\{(x,
y)\,|\,\,\frac{x^2}{a^2} + \frac{y^2}{b^2} =1,\,\,x \geq 0,\, y
\geq 0 \}$ \item let $\omega$ be a compact plane domain with the
smooth boundary, which lies in the (open) upper half-plane
${\mathbb R}_+^2$; let $\Omega$ be a torus in ${\mathbb R}^3$,
which appears as result of rotation of $\omega$ around the
$x$-axis. Then $\Omega_{L_0}\overset{\rm isom}=\omega$ and
$\partial \Omega_{L_0} \overset{\rm isom}=\partial \omega$.
\end{itemize}
\medskip

Actually, possible lack of simplicity is not an obstacle for
solving the problems IP 1--3 because their data determine not only
a copy of  $L_0$ but substantial additional information relevant
to the reconstruction of $\Omega$. Roughly speaking, the matter is
as follows. When we deal with each of these problems, the boundary
$\Gamma$ is given. By this, instead of the algebra ${\mathfrak N}[E_L,a]=
{\mathfrak N}_{L_0}$ generated by the family $a=\{{\rm clos\,} {\cal U}^t\}_{t \geq 0}$
of the sets reachable from {\it the whole} $\Gamma$ (see
(\ref{reachable_sets_U}) and (\ref{u^f=h-int})) we can invoke the
wider algebra ${\mathfrak N}[E_L,a'] \supset {\mathfrak N}[E_L,a]$ generated by the
much richer family $a'=\{{\rm clos\,} {\cal U}^t_\sigma\}_{t \geq 0, \,\sigma
\subset \Gamma} \supset a$ of the sets reachable from {\it any
patch} $\sigma \subset \Gamma$ of positive measure
\footnote{${\cal U}^t_\sigma$ consists of the solutions (waves) $u^f(t)$
produced by the boundary controls $f$ supported on $\sigma \times
[0, \infty)$}. As result, although the equality ${\mathfrak N}_{L_0}={\frak
L}_\infty(\Omega)$ may be violated by symmetries, the equality
${\mathfrak N}[E_L,a']={\frak L}_\infty(\Omega)$ always holds, whereas the
wave spectrum $\Omega_{{\mathfrak N}[E_L, a']}$ turns out to be isometric
to $\Omega$. The latter is the key fact, which makes the
reconstruction possible: see \cite{BSobolev} for detail.

\subsection{Comments and remarks}
{\bf A look at isospectrality}\,\,\, As at the end of sec 2.3, let
${\rm spec}\,L=\{\lambda_k\}_{k=1}^\infty$ be the spectrum of the
Dirichlet Laplacian on $\Omega$. The question: "Does ${\rm
spec}\,L$ determine $\Omega$ up to isometry?" is a version of the
classical M.Kac's drum problem \cite{Kac}. The negative answer is
well known (see, e.g., \cite{BCDS_isosp_dom}) but, as far as we
know, the satisfactory description of the set of isospectral
manifolds is not obtained yet. The following is some observations
in concern with such a description.

Assume that we deal with a simple $\Omega$. In accordance with
Theorem 1, such a manifold is determined by any unitary copy
$\widetilde L_0$ of the operator $L_0 \subset L$. If the spectrum
of $L$ is given, to get such a copy it suffices to possess the
Fourier image $\widetilde {\cal D} =U{\cal D}$ of the harmonic
subspace in $\widetilde {\cal H} = \emph{l}_2$: see
(\ref{tildeL0}) \footnote{It is the fact, which is exploited in
\cite{BKac89}}. In the mean time, as is evident, if $\Omega$ and
$\Omega'$ are isometric, then the corresponding images are
identical: $\widetilde {\cal D}=\widetilde {\cal D}'$. Therefore,
$\Omega$ and $\Omega'$ are isospectral but not isometric iff
$\widetilde {\cal D}\not=\widetilde {\cal D}'$. In other words,
the subspace $\widetilde {\cal D}$ is a relevant "index" that
differs the isospectral manifolds.

To be a candidate on the role of the harmonic functions image,
which is admissible for the given $\widetilde L={\rm
diag\,}\{\lambda_1, \lambda_2, \,\dots\}$ (see (\ref{diag_L})), a
subspace $\widetilde {\cal D} \subset \emph{l}_2$ has to possess
the following properties:
\begin{enumerate}\item a lineal set ${\cal L}_{\widetilde {\cal D}} :={\widetilde
L}^{-1}\left[\emph{l}_2 \ominus \widetilde {\cal D}\right]$ is
dense in $\emph{l}_2$ (see (\ref{tildeL0})), whereas replacement
of $\widetilde {\cal D}$ by any wider subspace ${\widetilde{\cal
D}}^\prime \supset \widetilde {\cal D}$ leads to the lack of
density: ${\rm clos\,} {\cal L}_{\widetilde {\cal
D}^\prime}\not=\emph{l}_{2}$ \item extending an operator
$\widetilde L|_{{\cal L}_{\widetilde {\cal D}}}$ by Friedrichs,
one gets $\widetilde L$ (see (\ref{tildeL0})).
\end{enumerate}
In the mean time, taking {\it any} subspace $\widetilde {\cal D}
\subset \emph{l}_2$ provided 1 and 2 \footnote{such subspaces do
exist (M.M.Faddeev, private communication)}, one can construct a
symmetric operator $\widetilde L_0$ by (\ref{tildeL0}) and then
find its wave spectrum $\Omega_{\widetilde L_0}$ as a candidate to
be a drum. However, the open question is whether such a "drum" is
human (is a manifold).
\bigskip

\noindent{\bf Model}\,\,\,Once again, let $\Omega$ be simple. By
Lemma \ref{NL_0=L_infty}, the algebra ${\mathfrak N}_{L_0}$ is
{\it cyclic}, i.e., possesses the elements $g \in {\cal
H}=L_2(\Omega)$ such that
$${\rm clos\,} \left\{Pg\,|\,\,P \in {\rm
Proj\,}{\mathfrak N}_{L_0}\right\}\,=\,{\cal H}\,.$$ This enables one to realize
elements of ${\cal H}$ as functions on the wave spectrum of ${\mathfrak N}_{L_0}$
by the following scheme:
\begin{itemize}\item Fix a cyclic $g \in {\cal H}$ and endow
$\Omega_{L_0}$ with a measure $\mu$ as follows. For a maximal
eikonal $\tau=\int_0^\infty t\,dP^t_\tau \in \Omega_{L_0}$ and a
ball $B_r[\tau]:=\{\tau' \in
\Omega_{L_0}\,|\,\,\|\tau'-\tau\|<r\}$ put
\begin{equation}\label{measure mu def}
\mu\left(B_r[\tau]\right):=\left(P_\tau^r g, g\right)_{\cal H}
\end{equation} and
then extend $\mu$ to the Borel subsets of $\Omega_{L_0}$. As is
easy to check, the equality
\begin{equation}\label{measure mu representation}
\mu\left(B_r[\tau]\right)=\int_{\Omega^r[x]}|g|^2\,\,d\,\rm vol
\end{equation}
holds, where $\tau$ and $x$ are related through the isometry
$i:\Omega_{L_0} \to \Omega$ established by Theorem
\ref{Theorem_1}: $x=i(\tau)$.

So, we have a {\it model space} ${\cal H}_{\rm mod}:=L_{2,\,
\mu}(\Omega_{L_0})$. \item For a $y \in {\cal H}$, define its
image $I y \in {\cal H}_{\rm mod}$ by
\begin{equation}\label{Image def}
\left(I y\right)(\tau):=\underset{r \to 0} \lim\,
\frac{\left(P_\tau^r y, g\right)_{\cal H}}{\left(P_\tau^r g,
g\right)_{\cal H}}\,.
\end{equation}
The relations (\ref{measure mu def}) and (\ref{measure mu
representation}) easily imply
\begin{equation}\label{Image vis integrals}
\left(I y\right)(\tau)=\underset{r \to 0} \lim\,
\frac{\int_{\Omega^r[x]} y\,\bar g\,\,d\,\rm
vol}{\int_{\Omega^r[x]} |g|^2\,\,d\,\rm
vol}\,=\,\left(g^{-1}y\right)(x)\,
\end{equation}
and, hence, the image map $I$ is a unitary operator from ${\cal
H}$ onto ${\cal H}_{\rm mod}$. \end{itemize} An operator $$L^{\rm
mod}_0:=I L_0 I^*$$ can be regarded as a functional model of the
operator $L_0$ on its wave spectrum that we call a {\it wave
model}. In the case under consideration, we have $L^{\rm
mod}_0=g^{-1}L_0 g$, i.e., this model is just a gauge transform of
the original. As such, it is a {\it local} operator:
\begin{equation}\label{locality}
{\rm supp\,} L^{\rm mod}_0 w \subset {\rm supp\,} w
\end{equation}
holds for $w \in {\rm Dom\,} L^{\rm mod}_0$.
\bigskip

\noindent{\bf Conjectures}\,\,\,

{\bf 1.}\,\,\,We suggest and hope that ${\cal H}_{\rm mod}$ and
$L^{\rm mod}_0$ do exist for a wide class of symmetric
semi-bounded operators $L_0$, the locality property
(\ref{locality}) being held. In contrast to the known models (see,
e.g., \cite{Shtraus}), this one has good chances to be of the real
use for applications.

The principal point (and difficulty) is to attach the invariant
meaning to the limit (\ref{Image def}). Presumably, it can be done
in the framework of the representation
\begin{equation}\label{N=integral N}
{\mathfrak N}_{L_0}\,=\,\oplus \int_{\Omega_{L_0}}{{\mathfrak N}}_\tau\,d\mu(\tau)
\end{equation}
in the form of a 1-st kind von Neumann algebra. Such a
representation diagonalizes ${\mathfrak N}_{L_0}$ and is expected
to be valid for $L_0$'s coming from mathematical physics. As an
encouraging example, the Maxwell system in electrodynamics can be
mentioned: see \cite{BIP07}. Note that Maxwell's ${\mathfrak
N}_{L_0}$ is noncommutative.

{\bf 2.}\,\,\,A question of independent interest is whether any
von Neumann algebra with space extension (in particular,
${\mathfrak N}_{L_0}$) is of the form (\ref{N=integral N}). Also, in
addition to the metric $\|\tau-\tau'\|$, it is reasonable to look
for more subtle structures on $\Omega_{L_0}$ like tangent spaces,
differentiable structure, etc.

{\bf 3.}\,\,\,One more attractive option is to construct a {\it
wave model} of the abstract Green system satisfying Ryzhov's
axiomatics, by the scheme
$$\{ {\cal H},{\cal G};
A,\Gamma_0, \Gamma_1\} \Longrightarrow L_0 \Longrightarrow
\Omega_{L_0} \Longrightarrow \{ {\cal H}_{\rm mod}, {\cal G};
A^{\rm mod},\Gamma_0^{\rm mod}, \Gamma_1^{\rm mod}\}$$ with ${\cal
H}_{\rm mod}=\oplus \int_{\Omega_{L_0}}{\cal H}_\tau\,d \mu(\tau)$
and a {\it local} $ A^{\rm mod}$. An intriguing point is that the
boundary $\partial \Omega_{L_0}$ is well defined, so that there is
a chance to realize $\Gamma_{0, 1}^{\rm mod}$ as the "true" trace
operators. Such a model would provide a canonical realization of
the original system, the realization being determined by its Weyl
function $M$ and, hence, relevant to inverse problems.
\bigskip

\noindent{\bf Unbounded case}\,\,\,Return to the sec 1.1. If the
assumption (\ref{C<infty}) is cancelled, the set ${\rm
Eik\,}\mathfrak N [E,a]$ is still well defined but unbounded. In
particular, the case that {\it all} eikonals $\tau=\int_0^\infty
t\,dP^t$ are unbounded operators is realized in applications: for
instance, it holds if we deal with a noncompact simple manifold
$\Omega$. As result, even though in the mentioned example the
relevant maximal eikonals do exist, in the general situation we
have to correct the definition of the wave spectrum
$\Omega_{\mathfrak N [E,a]}$. A possible way out is to deal with
the {\it regularized} eikonals $\tau=\int_0^\infty
\frac{t}{1+\alpha t}\,dP^t$ with a fixed $\alpha>0$ and thus
reduce the situation to the bounded case.
\bigskip

\noindent{\bf A bit of philosophy}\,\,\,In applications, the
external observer pursues the goal to recover a manifold $\Omega$
via measurements at the boundary $\Gamma$. The observer prospects
$\Omega$ with waves $u^f$ produced by boundary controls. These
waves propagate into the manifold, interact with its inner
structure and accumulate information about the latter. The result
of interaction is also recorded at $\Gamma$. The observer has to
extract the information from the recorded.

By the rule of game in IPs, the manifold itself is unreachable in
principle. Therefore, the only thing the observer can hope for, is
to construct from the measurements an {\it image} of $\Omega$
possibly resembling the original. By the same rule, the only
admissible material for constructing is the waves $u^f$.

To be properly formalized, such a look at the problem needs two
things:
\begin{itemize}
\item an object that codes exhausting information about $\Omega$
and, in the mean time, is determined by the measurements \item a
mechanism that decodes this information. \end{itemize} Resuming
our paper, the first is the minimal Laplacian $L_0$, whereas to
decode information is to determine its wave spectrum constructed
from the waves $u^f$. It is $\Omega_{L_0}$, which is a relevant
image of $\Omega$.
\medskip

The current paper develops an algebraic trend in the BC-method
\cite{BSobolev}, by which {\it to solve IPs is to find spectra of
relevant algebras}. An attempt to use this philosophy for solving
new problems would be quite reasonable. An encouraging fact is
that in all above-mentioned unsolved IPs of anisotropic elasticity
and electrodynamics, graphs with cycles, etc, the wave spectrum
$\Omega_{L_0}$ does exist. However, to recognize how it looks like
and verify (if true!) that $\Omega_{L_0} \overset{\rm isom}=
\Omega$  is difficult in view of very complicated structure of the
corresponding reachable sets ${\cal U}^t$.

\end{document}